\documentclass[12pt]{amsart}
 \usepackage{amsmath,amssymb}
\usepackage{graphicx}
 {\theoremstyle{definition}
 \newtheorem{remark}{Remark}
 
 }

 \newcommand{\gggg}{\gamma}

 \renewcommand{\bf}{\bfseries}
\newcommand{\I}{{\mathbf{1\!\!\!1}}}
\newcommand{\oo}{\infty}
\newcommand{\defeq}{\,{\buildrel {\rm def}\over =}}

 \newcommand{\aaa}{\alpha}

 \newtheorem{theorem}{Theorem}

 \renewcommand{\epsilon}{\varepsilon}

\newcommand{\tCC}{\widetilde{C}}
\newcommand{\tff}{\widetilde{f}}
\newcommand{\hCC}{\widehat{C}}
 
 \newcommand{\dis}{\displaystyle}
 \newcommand{\eee}{\epsilon}
 \newcommand{\mmm}{\mu}
 
\begin{document}
\title[the tail of the bilinear Hardy-Littlewood function]{A maximal inequality for the tail of the bilinear Hardy-Littlewood function}
 \author{I. Assani(*) and Z. Buczolich(**)}
\thanks{The first author acknowledges support by NSF grant DMS 0456627.
The second listed author was partially supported by NKTH and by the
Hungarian National Foundation for Scientific Research T049727.
\newline\indent \hfill\includegraphics{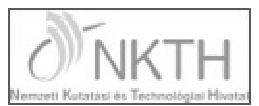}
}
 \begin{abstract}
 Let $(X,\mathcal{B}, \mu, T)$ be an ergodic dynamical system on a
 non-atomic finite measure space. We assume without loss of
 generality that $\mu(X)=1.$
 Consider the maximal function
 $\dis R^*:(f, g) \in L^p\times L^q \rightarrow R^*(f, g)(x) =
 \sup_{n\geq 1} \frac{f(T^nx)g(T^{2n}x)}{n}.$
 We obtain the following maximal inequality. For each $1<p\leq \infty$
 there exists a finite constant $C_p$ such that
 for each $\lambda >0,$ and
nonnegative functions $f\in L^p$ and $g\in L^1$
 $$\mu\{x: R^*(f,g)(x)>\lambda\}\leq
 C_p \bigg(
 \frac{\|f\|_p\|g\|_1}{\lambda}\bigg)^{1/2}.$$

 We also show that for each $\alpha>2$
the maximal function $R^*(f,g)$ is a.e. finite
 for pairs of functions $(f,g)\in
(L(\log L)^{2\alpha}, L^1)$.
 \end{abstract}
\maketitle
\section{Introduction}
Let $(X,\mathcal{B}, \mu, T)$ be an ergodic dynamical system on a
 non-atomic finite measure space. We assume without loss of
 generality that $\mu(X)=1.$

 In \cite{[AB]} we proved the following maximal inequality
about the maximal function
$\dis R^*(f, g)(x) =
 \sup_{n\geq 1} \frac{f(T^nx)g(T^{2n}x)}{n}.$ For each $1<p \leq \infty,$ there exists
 a finite constant $C_{p}'$ such that for each
 $\lambda>0,$ for every
$f\in L^p, f>1$ and $g\in L^1, g>1$

 \begin{equation}\label{meqq1}
 \mu\{x: R^*(f,g)(x) >\lambda\} \leq C_p' \bigg(
 \frac{\|f\|_p^p\|g\|_1}{\lambda}\bigg)^{1/2}.
 \end{equation}
 Furthermore the constant $C_p'$ behaves like $\frac{1}{p-1}$ when $p$
 tends to $1$. To be more precise, we will use that
there exists $\tCC'$ such that for any $1<p<2$ we have
\begin{equation}\label{cest1}
C_{p}'\leq \frac{\tCC'}{p-1}.
\end{equation}
Inequality \eqref{meqq1}
was enough to prove the a.e. convergence to zero of
 the tail $\dis \frac{f(T^nx)g(T^{2n}x)}{n}$ of the double recurrence averages $\dis
 \frac{1}{n}\sum_{k=1}^n
 f(T^kx)g(T^{2k}x)$ for pairs of functions $(f,g)$ in $L^p\times
 L^1$ (or $L^1\times L^p$) as soon as $p>1$. On the other hand, in \cite{[AB11]} the
 tail is used to show that these averages do not converge a.e. for
 pairs of $(L^1, L^1)$ functions.

 \vskip1ex

 During the 2007 Ergodic Theory workshop at UNC-Chapel Hill,
J.P. Conze asked if
 this inequality could be made homogeneous with respect to $f$ and $g$.
In this paper first we derive from \eqref{meqq1} the following homogeneous
 version.

 \begin{theorem}\label{th1}
For each $1<p<\infty$ there exists a finite constant $C_p$ such
that for each $f, g \geq 0$ and for all $\lambda> 0$ we have
\begin{equation}\label{meqq2}
\mu\big\{x: \sup_{n\geq 1} \frac{f(T^nx)g(T^{2n}x)}{n}
>\lambda\big\} \leq C_p \bigg(\frac{\|f\|_p\|g\|_1}{ \lambda
}\bigg)^{1/2},
\end{equation}
and there exists $\tCC$ such that for any
$1<p<2$ we have
\begin{equation}\label{cest2}
C_{p}\leq \frac{\tCC}{p-1}.
\end{equation}
 \end{theorem}

 At the same meeting a question
was raised about the a.e. finiteness
 of $R^*(f,g)$ for pairs of functions in $(L\log L, L^1).$
 Our second result is based on an adaptation of Zygmund's
 extrapolation method \cite{[Z]} (vol. II, ch. XII,
 pp. 119-120) to $R^*(f,g).$ With somewhat
 crude estimates we prove the following theorem.
\begin{theorem}\label{th2}
If $\alpha>2$ and the pair of nonnegative functions
$(f,g)$ belongs to\\ $ (L(\log L)^{2\alpha}, L^1)$
then
 $R^*(f,g)=
 \sup_{n\geq 1} \frac{f(T^nx)g(T^{2n}x)}{n}$ is a.e. finite.
\end{theorem}

\section{Proofs}

\begin{proof}[Proof of Theorem \ref{th1}]
First we can notice that the original inequality (1) is homogeneous
with respect to the $L^1$ function $g$. Indeed, a simple change of
variables shows that the case $g>t$ can easily be obtained from the
case $g>1$ with the same constant
$C_{p}'$. So
by approximating $g$ with $g_{n}(x)=\max \{ g(x),1/n \}$
we can see that
\eqref{meqq1} holds if the assumption $g> 1$ is replaced
by $g\geq 0$.
Without loss of generality we can also suppose
in the sequel that
$\|g\|_1 =1$.

\vskip1ex
If $\| f\| _{p}=0$ we have nothing to prove.
Otherwise, if we can show that
 \eqref{meqq2} holds for $\tff =f/\| f\| _{p}$ for all $\lambda >0$,
then this implies that it is true for $f$ as well for all $\lambda >0$.
Thus, we just need to prove \eqref{meqq2}
 for $f\in L^p$ with $\|f\|_p =1.$

Set
$$M=\mu\big\{x: \sup_{n\geq 1} \frac{f(T^nx)g(T^{2n}x)}{n}
>\lambda\big\}$$
and $h = \max\{f, 1\}.$ By
our remark about the assumption
$g\geq 0$ the maximal inequality \eqref{meqq1}
is applicable and
 we obtain that
 $M\leq C_p' \big(\frac{\|h\|_p^{p} }{\lambda}\big)^{1/2}$,
and \eqref{cest1} also holds for $1<p<2$.
 As $\|h\|_p\leq \|\I \|_p + \|f\|_p=2$ we have the estimate
$$M\leq 2^{p/2}  C_p' \bigg(\frac{1}{\lambda}\bigg)^{1/2}=
2^{p/2}  C_{p}'\bigg(\frac{\| f\| _{p}\| g\|
_{1}}{\lambda}\bigg)^{1/2},$$ with $C_{p}'$ satisfying \eqref{cest1}
for $1<p<2.$
 Therefore,  we obtain
$$\mu\big\{x: \sup_n \frac{f(T^nx)g(T^{2n}x)}{n}
>\lambda\big\}\leq
 2^{p/2}  C_{p}'
\bigg(\frac{\| f\| _{p}\| g\| _{1}}{\lambda}\bigg)^{1/2}\leq
%$$
%$$
C_{p}
\bigg(\frac{\| f\| _{p}\| g\| _{1}}{\lambda}\bigg)^{1/2}$$ with $C_{p}=
2^{p/2}  C_{p}'$ and from \eqref{cest1} it follows
that there exists $\tCC$ such that \eqref{cest2} holds for $1<p<2$.
\end{proof}

\vskip1ex

\begin{proof}[Proof of Theorem \ref{th2}]

 The starting point is \eqref{meqq2} and \eqref{cest2}.

 There exists a finite constant $\tCC$ such
that
 for every $1<p<2$,
for each $f, g \geq 0$ and for all $\lambda> 0$ we have
\begin{equation}\label{meqq3}
\mu\big\{x: \sup_n \frac{f(T^nx)g(T^{2n}x)}{n} >\lambda\big\} \leq
\frac{\tCC}{p-1} \bigg(\frac{\|f\|_p\|g\|_1}{\lambda }\bigg)^{1/2}.
\end{equation}

We can again assume without loss of generality that $\|g\|_1=1.$ We fix
the function $g$ and denote by $R^{*}(f)(x)$ the maximal function
$\sup_n \frac{f(T^nx)g(T^{2n}x)}{n}$.
Now we can rewrite \eqref{meqq3} as
\begin{equation}\label{*3}
\mu\big\{x: R^{*}(f)(x) >\lambda\big\} \leq
\frac{\tCC}{p-1} \bigg(\frac{\|f\|_p }{\lambda }\bigg)^{1/2}.
\end{equation}

 The important element for the extrapolation is the factor
$\frac{1}{p-1}$ in the above inequality.

Our goal is to prove that for
$\alpha>2$ there is $C_{\aaa}$ such that for any
 $f\in L(\log L)^{2\alpha}$ we have for
each $\lambda>0$
\begin{equation}\label{meqq4}
\mu\big\{x: R^{*}(f)(x)>\lambda\big\} \leq C_{\aaa}
\frac{1+ \big( \int |f|(\log^{+}
|f|)^{2\aaa}\big)^{1/2}}{{\lambda}^{1/2}} .
\end{equation}

Let $\gggg_j$ be a positive sequence of numbers such that $\dis
\sum_{j=0}^{\infty} \gggg_j = 1.$

The function $f$ being in $L(\log L)^{2\aaa}$ we have $\dis
\sum_{j=1}^{\infty} j^{2\aaa}2^j\mu\big\{x: 2^j\leq
f<2^{j+1}\big\}<\infty.$ We denote by $t_j $ the quantity $\mu\big\{
2^j\leq f<2^{j+1}\big\},$ by $f_j$ the function $2^j{\I}_{\big\{x:
2^j\leq f<2^{j+1}\big\}}$ and by $p_j$ the number $1+ \frac{1}{j}.$
We set $f_{0}(x)=f(x)$ if $0\leq f(x)<2$, otherwise we put
$f_{0}(x)=0$.
Then
\begin{equation}\label{*5x*}
f\leq 2\sum_{j=0}^{\oo}f_{j}.\end{equation}
We also have
\begin{equation}\label{foest}
\mu\big\{x: R^{*}(f_{0})(x)>\frac{\lambda\gggg_{0}}
{2}\big\} \leq
\end{equation}
$$\mu\big\{x: R^{*}(2\cdot \I_{X})(x)>\frac{\lambda\gggg_{0}}
{2}\big\}\leq \frac{4\| g\| _{1}}{\lambda\gggg_{0}}=
\frac{4}{\lambda \gggg_{0}}$$ by the standard maximal inequality for
the ergodic averages (see \cite{[K]} for instance).

For $j\geq 1$ by \eqref{*3} used with $p_{j}=1+\frac{1}j$
we obtain
\begin{equation}\label{fjest}
\mu\big\{x: R^{*}(f_{j})(x)>\frac{\lambda\gggg_{j}}
{2}\big\} \leq
\end{equation}
$$\tCC \frac{1}{(1+(1/j))-1}
\left (\frac{2^{j/2}[t_{j}]^{1/2p_{j}}}{(\lambda \gggg_{j}/2)^{1/2}}
 \right)\leq
 %$$
%$$
\sqrt 2\tCC
\frac{j2^{j/2}[t_{j}]^{1/2p_{j}}}{(\lambda \gggg_{j})^{1/2}}
.$$

 We choose
$\gggg_{0}=1/2$ and $\gggg_j = {\frac{C_{\gggg}}{j(\log(j+1))^{\beta}}}$ with
 $\beta>1$ and $C_{\gggg}$ such that $\dis \sum_{j=0}^{\infty} \gggg_j=1.$

Set $\dis \hCC=\frac{\sqrt 2\tCC}{C_{\gggg}^{1/2}}
.$

Using \eqref{*5x*} and
 adding \eqref{foest} and
\eqref{fjest}
for all $j$ we obtain
\begin{equation}\label{*0804}
\mu\big\{x: R^{*}(f)(x)>\lambda\big\} \leq \sum_{j=0}^{\infty}
\mu\big\{R^{*}(f_j)>\frac{\lambda\gggg_j}2\big\} \leq
 \frac{8}{\lambda }+ \sqrt 2\tCC
\sum_{j=1}^{\infty}\frac{j2^{j/2}\big[t_j\big]^{1/2p_j}}{
 {(\lambda\gggg_j)}^{1/2}} \leq
 \end{equation}

% \[\begin{aligned}
% &\mu\big\{x: R^{*}(f)(x)>\lambda\big\}\\
% &
$$\frac{8}{\lambda }+\hCC \frac{\sum_{j=1}^{\infty}j^{3/2}[\log
 (j+1)]^{\beta/2}2^{j/2}\big[t_j\big]^{1/2p_j}}
 {{\lambda}^{1/2}} =
%$$
%\\
% \end{aligned}
% \]
%$$
\frac{8}{\lambda }+\hCC\frac{A_{1}}{\lambda^{1/2} }.$$
To estimate $A_{1}$
denote by $J_{1}$ the set of those $j$ for which
$t_j^{1/2p_j}\leq 3^{-j}$.
Then
\begin{equation}\label{*3*}
\sum_{j\in J_{1}}j^{3/2}[\log
 (j+1)]^{\beta/2}2^{j/2}\big[t_j\big]^{1/2p_j}
\leq
\sum_{j=1}^{\oo}j^{3/2}[\log
 (j+1)]^{\beta/2}2^{j/2}3^{-j}\defeq C_{s}.
\end{equation}
If $j\not\in J_{1}$ then $t_{j}^{1/2p_{j}}>3^{-j}$,
that is, $$3>t_{j}^{\frac{-\frac{1}j}{2p_{j}}}=t_{j}^{\frac{1-(1+\frac{1}j)}
{2p_{j}}}=
t_{j}^{\frac{1}{2p_{j}}-\frac1{2}},$$
which implies
$t_j^{1/2p_j}< 3t_j^{1/2}.$
Hence
\begin{equation}\label{*3**}
\sum_{j\not\in J_{1}}j^{3/2}[\log
 (j+1)]^{\beta/2}2^{j/2}\big[t_j\big]^{1/2p_j}
\leq
3\sum_{j=1}^{\oo}j^{3/2}[\log
 (j+1)]^{\beta/2}2^{j/2}\big[t_j\big]^{1/2} \defeq B_{1}.
\end{equation}
Suppose that $\aaa>\delta >2$.
By rewriting and applying the Cauchy--Schwartz inequality
we obtain with a suitable constant $C_{\delta }$ that
 $$B_{1}=3\sum_{j=1}^{\infty}\big[j^{3/2}
 {j^{-\delta}}\big] j^{\delta}\big[\log
 (j+1)]^{\beta/2}2^{j/2}\big[t_j\big]^{1/2}\leq$$
 $$3\big[\sum_{j=1}^{\infty}j^{3-2\delta}\big]^{1/2}\big[\sum_{j=1}^{\infty}j^{2\delta}\big[\log
 (j+1)\big]^{\beta}2^jt_j\big]^{1/2}=$$
$$C_{\delta }
\big[\sum_{j=1}^{\infty}j^{2\delta}\big[\log
 (j+1)\big]^{\beta}2^jt_j\big]^{1/2}\defeq B_{2}.$$
There exists $C_{\delta ,\aaa,\beta}$
such that for all $j=1,2,...$
$$\big[\log
 (j+1)\big]^{\beta}
\leq C_{\delta ,\aaa,\beta} j^{2(\aaa-\delta )}.$$
Hence,
\begin{equation}\label{12**}
B_{1}\leq B_2 \leq C_{\delta }C_{\delta ,\aaa,\beta}
\left (\int |f|(\log^{+} |f|)^{2\aaa}d\mu\right)^{1/2}.
\end{equation}

By
(\ref{*0804}-\ref{12**})
we have
$$
\mu \{ x: R^{*}(f)(x)>\lambda \}
\leq
%$$
%$$
\hCC \frac{C_{s}+ C_{\delta } C_{\delta,\aaa,\beta }(\int |f|
(\log^{+}|f|)^{2\aaa}
d\mu)^{1/2}
}
{\lambda ^{1/2}}
$$
this implies \eqref{meqq4} with a suitable $C_{\aaa}$.

\end{proof}

\begin{remark}
Inequality \eqref{meqq4} implies also that for the pair of
 nonnegative functions $(f, g)$ in $(L(\log L)^{2\alpha}, L^1)$ we
 have
 \begin{equation}\label{*78}
 \lim_n \frac{f(T^nx)g(T^{2n}x)}{n}= 0.
  \end{equation}

Indeed, consider
 a sequence of bounded functions $0\leq f_M\leq f$ converging monotone
increasingly to $f\in L(\log
 L)^{2\alpha}.$ Then we have
 \begin{equation}\label{*7*}
 \dis \lim_n
 \frac{f_M(T^nx)g(T^{2n}x)}{n}=0.
 \end{equation}
 Given $\eee\in (0,1)$ choose $M$ so large that
\begin{equation}\label{**7*}
I(M,\eee,1/2)\defeq\left (\int \frac{2}{\eee^{2}}|f-f_{M}|(\log^{+}\frac{2}{\eee^{2}}|f-f_{M}|)^{2\aaa}
d\mu \right )^{1/2}<1.
\end{equation}
Then
$$\mmm\{ x:\limsup_{n\to\oo}\frac{f(T^{n}x)g(T^{2n}x)}n>\eee  \}\leq$$
$$\mmm\{ x:\limsup_{n\to\oo}\frac{(f-f_{M})(T^{n}x)g(T^{2n}x)}n>\frac{\eee}2  \}
+\mmm\{ x:\limsup_{n\to\oo}\frac{f_{M}(T^{n}x)g(T^{2n}x)}n>\frac{\eee}{2}  \}\leq$$
(by using \eqref{*7*})
$$\mmm\{ x:R^{*}((f-f_{M}),g)(x)>\frac{\eee}2  \}
=\mmm\{ x:R^{*}(\frac{2}{\eee^{2}}(f-f_{M}),g)(x)>\frac 1{\eee}  \}\leq$$
(by using \eqref{meqq4} and \eqref{**7*})
$$C_{\aaa}\sqrt{\eee}(1+I(M,\eee,1/2))\leq 2C_{\aaa}\sqrt{\eee}.$$
Since this holds for any $\eee\in (0,1)$ we obtained \eqref{*78}.

\end{remark}

\vskip1ex

(*) Idris Assani - Department of Mathematics- University of North
Carolina at Chapel Hill-email: assani@email.unc.edu

\vskip1ex

(**) Zolt\'an Buczolich- Department of Analysis, E\"otv\"os Lor\'and
University, P\'azm\'any P\'eter S\'et\'any 1/c, 1117 Budapest,
Hungary -email: buczo@cs.elte.hu
\end{document}